# A counting method for finding rational approximates to arbitrary order roots of integers

Ashok Kumar Gupta, Department of Electronics and Communication,
Allahabad University, Allahabad - 211 002, India
(e-mail: ashok@mri.ernet.in)

Ashok Kumar Mittal, Department of Physics,
Allahabad University, Allahabad – 211 002, India
(e-mail:mittal@mri.ernet.in)

ABSTRACT:

It is shown that for finding rational approximates to $m^{th}$ root of any integer to any accuracy, one only needs the ability to count and to distinguish between $m$ different classes of objects. To every integer N can be associated a 'replacement rule' that generates a word **W\*** from another word **W** consisting of symbols belonging to a finite 'alphabet' of size $m$. This rule applied iteratively on almost any initial word **W₀**, yields a sequence of words **{Wᵢ}** such that the relative frequency of different symbols in the word **Wᵢ**, approaches powers of the $m^{th}$ root of N as $i$ tends to infinity.

Let us consider the problem of finding rational approximates to arbitrary order (irrational) roots of integers. For the square roots we have a well-known standard method, but we do not have any comparable procedure for cube roots and higher order roots. They can only be calculated by making use of logarithms and anti-logarithms. We present in this paper a conceptually simple, but remarkable, method for calculating arbitrarily close rational approximates to arbitrary order roots of integers. This method involves only the ability to count and to distinguish between $m$ different classes of objects, if the $m^{th}$ order root is desired.

We shall begin with a simple example of calculating the square root of 2. We begin with a binary alphabet {0,1}. Let the symbols '0' and '1' be replaced by the following replacement rules.

$$0 \rightarrow 01$$
$$1 \rightarrow 100 \tag{1}$$

We begin with the symbol '0' as the seed and continue to apply the above replacement rules repeatedly on each symbol in the sequence. We generate the following (aperiodic) sequence of sequences:

|    |                                    | n0 | n1 | n0/n10   |
|----|------------------------------------|----|----|----------|
| 0. | 0                                  | 1  | 0  | infinity |
| 1. | 01                                 | 1  | 1  | 1        |
| 2. | 01100                              | 3  | 2  | 1.5      |
| 3. | 011001000101                       | 7  | 5  | 1.4      |
| 4. | 0110010001100010101100 1100         | 17 | 12 | 1.4165   |

Here n0 and n1 are the numbers of 0's and 1's in the sequence. If we continue to generate these sequences, we get the following series of rational approximates to the square root of 2:

|    | n0     | n1     | n0/n1            |
|----|--------|--------|------------------|
| 1  | 1      | 1      | 1                |
| 2  | 3      | 2      | 1.50000000000000 |
| 3  | 7      | 5      | 1.40000000000000 |
| 4  | 17     | 12     | 1.41666666666667 |
| 5  | 41     | 29     | 1.41379310344828 |
| 6  | 99     | 70     | 1.41428571428571 |
| 7  | 239    | 169    | 1.41420118343195 |
| 8  | 577    | 408    | 1.41421568627451 |
| 9  | 1393   | 985    | 1.41421319796954 |
| 10 | 3363   | 2378   | 1.41421362489487 |
| 11 | 8119   | 5741   | 1.41421355164605 |
| 12 | 19601  | 13860  | 1.41421356421356 |
| 13 | 47321  | 33461  | 1.41421356205732 |
| 14 | 114243 | 80782  | 1.41421356242727 |
| 15 | 275807 | 195025 | 1.41421356236380 |
| 16 | 665857 | 470832 | 1.41421356237469 |

These ratios are quadratically converging to $(2)^{1/2} = 1.41421356237310\ldots$ It should be remembered that we have used only replacement rules (1) and counting abilities to get arbitrarily close rational approximates to $(2)^{1/2}$.

In order to get $(3)^{1/2}$, we modify the replacement rules as follows:

    0 ➔ 01
    1 ➔ 1000                                                                       (2)

Again we proceed with the starting symbol '0', generate the sequences as before and count the number of 0's and 1's. The ratio n0/n1 converges to $(3)^{1/2} = 1.73205080756888\ldots$

We can get the square root of any positive integer N by the same method by using the replacement rules:

    0 ➔ 01
    1 ➔ 100…0               (1 followed by N 0's)                      (3)

The method can now be generalized to the cube roots. As an example take the cube root of 2. Let us now have three symbols, 0, 1 and 2, and the replacement rules as follows;

    0 ➔ 01
    1 ➔ 12                                                                                              (4)
    2 ➔ 200

Now the following aperiodic sequence of sequences is generated:

|   |   | n0 | n1 | n2 | n0/n1 | n1/n2 |
|---|---|----|----|----|-------|-------|
| 0. | 0 | 1 | 0 | 0 | infinity | NaN |
| 1. | 01 | 1 | 1 | 0 | 1 | infinity |
| 2. | 0112 | 1 | 2 | 1 | 0.5 | 2 |
| 3. | 011212200 | 3 | 3 | 3 | 1 | 1 |
| 4. | 011212200122002000101 | 9 | 6 | 6 | 1.5 | 1 |

If we continue to generate these sequences, we get the following series of rational approximates to the cube root of 2:

|    | n0    | n1    | n2    | n0/n1            | n1/n2            |
|----|-------|-------|-------|------------------|------------------|
| 1  | 1     | 1     | 0     | 1                | infinity         |
| 2  | 1     | 2     | 1     | 0.50000000000000 | 2                |
| 3  | 3     | 3     | 3     | 1                | 1                |
| 4  | 9     | 6     | 6     | 1.50000000000000 | 1                |
| 5  | 21    | 15    | 12    | 1.40000000000000 | 1.25000000000000 |
| 6  | 45    | 36    | 27    | 1.25000000000000 | 1.33333333333333 |
| 7  | 99    | 81    | 63    | 1.22222222222222 | 1.28571428571429 |
| 8  | 225   | 180   | 144   | 1.25000000000000 | 1.25000000000000 |
| 9  | 513   | 405   | 324   | 1.26666666666667 | 1.25000000000000 |
| 10 | 1161  | 918   | 729   | 1.26470588235294 | 1.25925925925926 |
| 11 | 2619  | 2079  | 1647  | 1.25974025974026 | 1.26229508196721 |
| 12 | 5913  | 4698  | 3726  | 1.25862068965517 | 1.26086956521739 |
| 13 | 13365 | 10611 | 8424  | 1.25954198473282 | 1.25961538461538 |
| 14 | 30213 | 23976 | 19035 | 1.26013513513514 | 1.25957446808511 |
| 15 | 68283 | 54189 | 43011 | 1.26008968609865 | 1.25988700564972 |

It may be observed that both of the ratios converge to cube-root of 2 = 1.25992..

The generalization to the cube root of any integer N is straightforward. The replacement rules in this case are:

> 0 → 01
> 1 → 12  (5)
> 2 → 200..0   (2 followed by N 0's)

The generalization to find the $m^{th}$ root of an integer N is now possible. The replacement rules are:

$0 \to 01$
$1 \to 12$
$2 \to 23$
$3 \to 34$            (6)
.
.
.
$(m-1) \to (m-1)\ 0\ 0\ 0\ ...0$     ($(m-1)$ followed by N 0's)

The above procedure for obtaining rational approximates to the $m^{th}$ root of integer N does not involve any mathematical operation more advanced than counting. The above procedure can be justified mathematically as follows:

Let $A = \{A_1, A_2, ...A_m\}$ be a finite set. Let $A^* = \Sigma_p A \times A \times .....\times A$ (p times). An element of $A^*$ is called a 'word' obtained from the 'alphabet' A.

Consider the replacement rule $R: A \to A^*$ given by

$$A_1 \to A_1 A_2$$
$$A_2 \to A_2 A_3$$
.
.          (7a)

$$A_{m-1} \to A_{m-1} A_m$$
$$A_m \to A_m A_1........A_1$$      ($A_1$ is repeated N times)

Let $W = A_{s1}A_{s2}.....A_{sq}$ where $s_k$ $(k=1,2,........,q) \in \{1,2,........,m\}$. Then $W \in A^*$. The replacement rule R induces a mapping $R^*: A^* \to A^*$ defined by

$$\begin{aligned} W^* &= R^*(W) \\ &= R^*(A_{s1}A_{s2}........A_{sq}) \\ &= R(A_{s1})R(A_{s2})........R(A_{sq}) \end{aligned}$$   (7b)

Let $W_0 \in A^*$. Then,

$$W_i = (R^*)^i (W_0) = R^*((R^*)^{i-1}(W_0)) = R^*(W_{i-1})$$     (8)

denotes the word obtained by i times repeated application of the replacement rule (7) on the initial word $W_0$.

Let **n**: $A^* \to V^m$ be a mapping, which assigns to a word in $A^*$, a vector

$$\mathbf{n}(W) = (n_1(W), n_2(W), \ldots, n_j(W), \ldots n_m(W))^T \qquad (9)$$

in an m-dimensional vector space $V^m$ such that $n_j(W)$ is a non-negative integer denoting the frequency of occurrence of $A_j$ in the word W.

The replacement rule (7) and the mapping **n** induces a m x m matrix **R** which maps vectors in $V^m$ into $V^m$ such that

$$\mathbf{n}(W^*) = \mathbf{n}(R^*(W)) = \mathbf{R}(\mathbf{n}(W)) \qquad (10)$$

It follows from (7) - (10) that

$$\mathbf{R} = \begin{bmatrix} 1 & 0 & 0 & \ldots & 0 & N \\ 1 & 1 & 0 & \ldots & 0 & 0 \\ 0 & 1 & 1 & \ldots & 0 & 0 \\ . & . & . & \ldots & . & . \\ . & . & . & \ldots & . & . \\ 0 & 0 & 0 & \ldots & 1 & 1 \end{bmatrix} \qquad (11)$$

Equations (9) and (11) imply

$$\mathbf{n}(W_i) = \mathbf{n}(R^*(W_{i-1})) = \mathbf{R}(\mathbf{n}(W_{i-1})) = \mathbf{R}^i(\mathbf{n}(W_0)) \qquad (12)$$

For almost any $W_0$, $\mathbf{n}(W_i)$ tends, as $i \to \infty$, to multiples of the eigenvector of **R** corresponding to the eigenvalue with maximum absolute value[1]. The largest eigenvalue of the matrix **R** is $1+ N^{1/m}$. The eigenvector corresponding to this eigenvalue is given by $[\lambda^{m-1}, \lambda^{m-2}, \ldots, \lambda, 1]^T$, where $\lambda = N^{1/m}$. Thus one finds that for almost any initial word $W_0$,

$$\lim_{i \to \infty} n_j(W_i)/n_{j+1}(W_i) = \lambda = N^{1/m} \qquad j = 1, 2, \ldots, m-1 \qquad (13)$$

Eqn (13) shows that rational approximates to arbitrary order root of any integer can be obtained to any accuracy only by operations consisting of replacement of symbols and counting symbols. This method does not require any 'advanced' operation such as multiplication, division, logarithms and anti-logarithms, although the justification of the method uses 'advanced' concepts like the theory of matrices. Of course, 'advanced operations' like multiplication, division, logarithms and anti-logarithms can also be reduced to long chains of counting processes. But the simplicity of the method given here, as compared with any standard method if broken down into operations of counting, appears to be quite non-trivial.

The numbers that can be obtained geometrically using only a ruler and a compass are called constructible numbers[2]. It is known that the set of constructible numbers is $Q(\sqrt{m_1}, \sqrt{m_2}, \sqrt{m_3}, \ldots, \sqrt{m_N})$, the extension field of rational numbers by square roots of integers. Higher roots such as cube roots etc, in general, are not constructible. The replacement and counting method of this paper shows how to construct these roots by a method more primitive than the use of ruler and compass. The set of numbers that may be constructed by this method, which may be called *constructible by replacement and counting,* contains the set of (geometrically) constructible numbers $Q(\sqrt{m_1}, \sqrt{m_2}, \ldots, \sqrt{m_N})$.